# Mathematical Insights in the Pioneering Educational Project FSTF

**Enzo Bonacci**
Department of Mathematics and Physics, Scientific High School "G.B. Grassi"
Latina (Italy)
E-mail: enzo.bonacci@liceograssilatina.org

**Abstract**

*The educational project "From Soccerene to Fullerene" was a pioneering interdisciplinary experience accomplished by six pupils who discussed it in the 2013 Italian school-leaving examination. The author illustrated the FSTF in the 2014 national conference of the Italian Physical Society (Pisa) eliciting a genuine interest. Here we highlight the FSTF's mathematical aspects.*

*Key Words: Plane and solid geometry, Archimedean solids, Deformable solids, Geometry and physics*

*MSC(2010): 97G40, 51P05, 74-02, 97M50*

**INTRODUCTION:**
During their pre-university year, six brilliant and motivated students joined the educational project "From Soccerene to Fullerene" (acronym FSTF) proposed by five teachers from the Departments of *Physical Education* and of *Mathematics and Physics* at the Scientific High School "G.B. Grassi" in Latina (Italy). The pupils investigated the *fullerene* buckyball ($C_{60}$) also known as *soccerene* for its football-shaped structure (truncated icosahedron). From plastic sheets they cut twelve pentagons and twenty hexagons, with a calculated edge, sewing them together on a polystyrene sphere as faces of a soccer ball; they repeated the same operation twice, making two almost identical soccer-like balls. That multidisciplinary experience started on March 2013 and ended on June 2013. It was satisfactorily discussed in the 2013 school-leaving examination and raised a great interest in the 100[th] national conference of the Italian Physical Society in Pisa, where it was presented the following year [1], [5]. Now we briefly illustrate the FSTS, focusing on its mathematical features.

**PREMISES OF THE FSTF:**
The teachers submitted the following four questions, asking the students to seek for the answers by themselves:
  (1) Why do we play soccer with a sphere?
  (2) Why is the soccer ball a truncated icosahedron?
  (3) What is the soccerene or fullerene?
  (4) Why is the fullerene so important?

**Why do we play soccer with a sphere?**
The pupils verified that, as fully expounded by Ludwig and Guerrerio [14], the physics of rubber elastic bodies [20] requires a perfectly spherical ball in order to avoid unpredictable trajectories after the rebound. In this way, the trajectory after the impact is independent of the point of contact between the ball and the hit surface.

**Why is the soccer ball a truncated icosahedron?**
The students detected that the curricular *Platonic solids* (tetrahedron, cube, octahedron, dodecahedron, icosahedron) were not suitable for constructing a spherical ball. They discovered how Archimedes introduced thirteen semi-regular convex polyhedra whose faces are two or more types of regular polygons. Then the pupils recognized that, as clarified by Ludwig and Guerrerio [14], only two *Archimedean solids* are closely approximated by a sphere:
(1) Rhombicosidodecahedrons: with volume 94% of the circumscribed sphere (12 pentagons, 30 squares and 20 triangles).
(2) Truncated icosahedrons: with volume 87% of the circumscribed sphere (12 pentagons, 20 hexagons).

The students understood that, between them, the best choice for building a soccer ball was the *truncated icosahedron*, with less edges (90 Vs. 120) by the Descartes-Euler polyhedral formula stating that the number of edges is equal to the number of faces plus the number of vertices minus two ($E = F + V - 2$) for any convex polyhedron [6], [21].

**What is the soccerene or fullerene?**
The pupils studied that, for its hexagonal texture interspersed by pentagonal asymmetries, the soccer ball is a miniaturized version of the geodesic domes patented in 1954 by the architect *Richard Buckminster Fuller* [15]. They learned that, to emphasize the structural analogy, the molecule of 60 carbon atoms was named *buckminsterfullerene* [13] and that *soccerene* [12] is another definition of the same $C_{60}$ molecule (IUPAC name *($C_{60}$-$I_h$)[5,6]fullerene*) whose shape resembles a soccer ball; hence the title of our project [19]. We remark that the pupils' team was already accustomed with the hexagonal geodesic dome of the planetarium "Livio Gratton" ([2], [3], [4]) located in the Scientific High School "G.B. Grassi" of Latina (Figs. 1, 2 and 3).

**Why is the fullerene so important?**
The students comprehended that the fullerene is a topical scientific issue, as explained by Porter [18] and remarked by Giuliano [9, 11]. They found that also the related new materials and nanotechnologies [16] are surely worthy of investigation in a scientific high school.

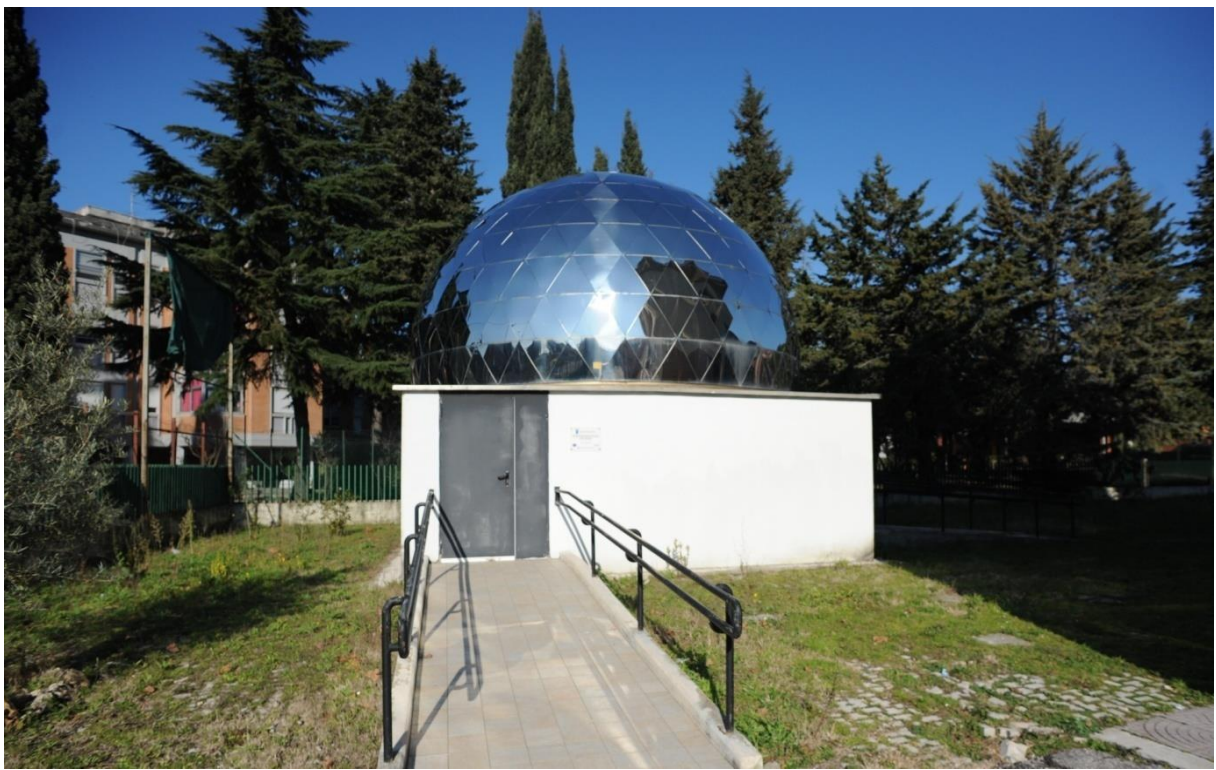

**Figure 1.** The planetarium "Livio Gratton" at the Scientific High School "G.B. Grassi" in Latina (Italy).

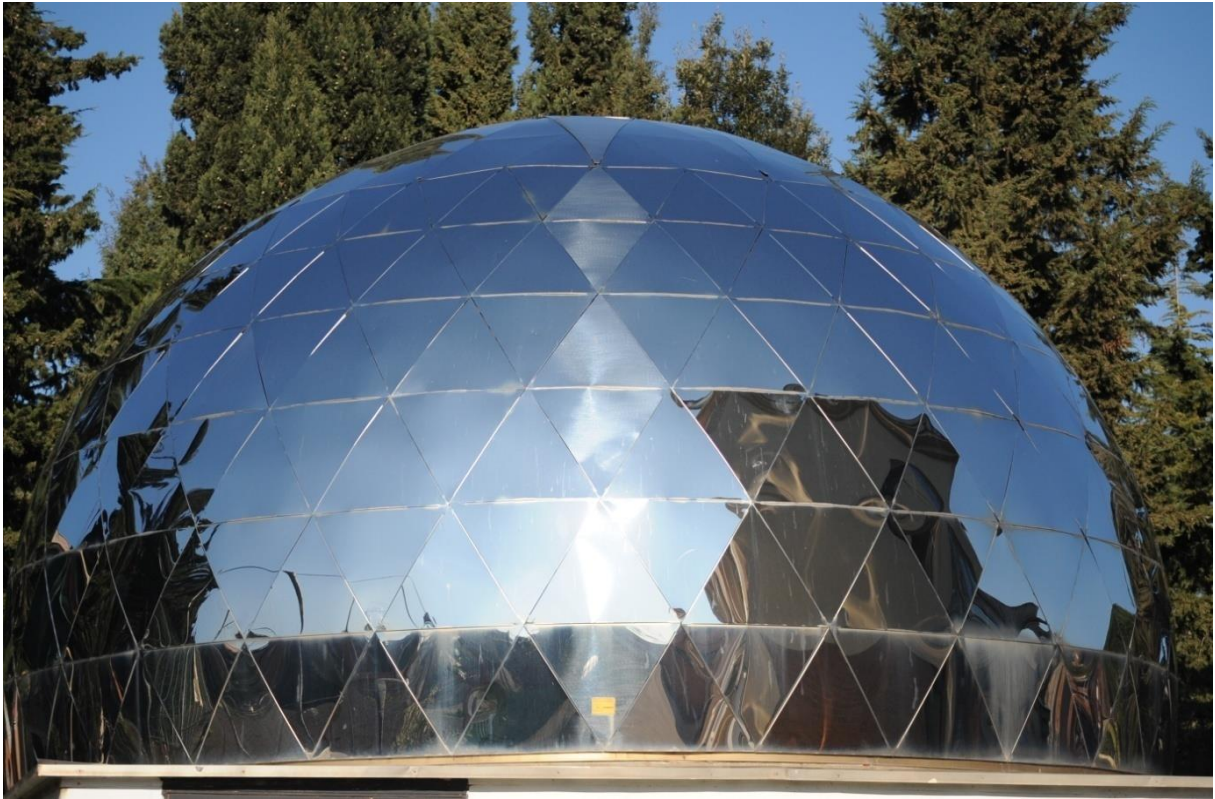

**Figure 2.** The geodesic dome of the planetarium "Livio Gratton" at the Scientific High School "G.B. Grassi".

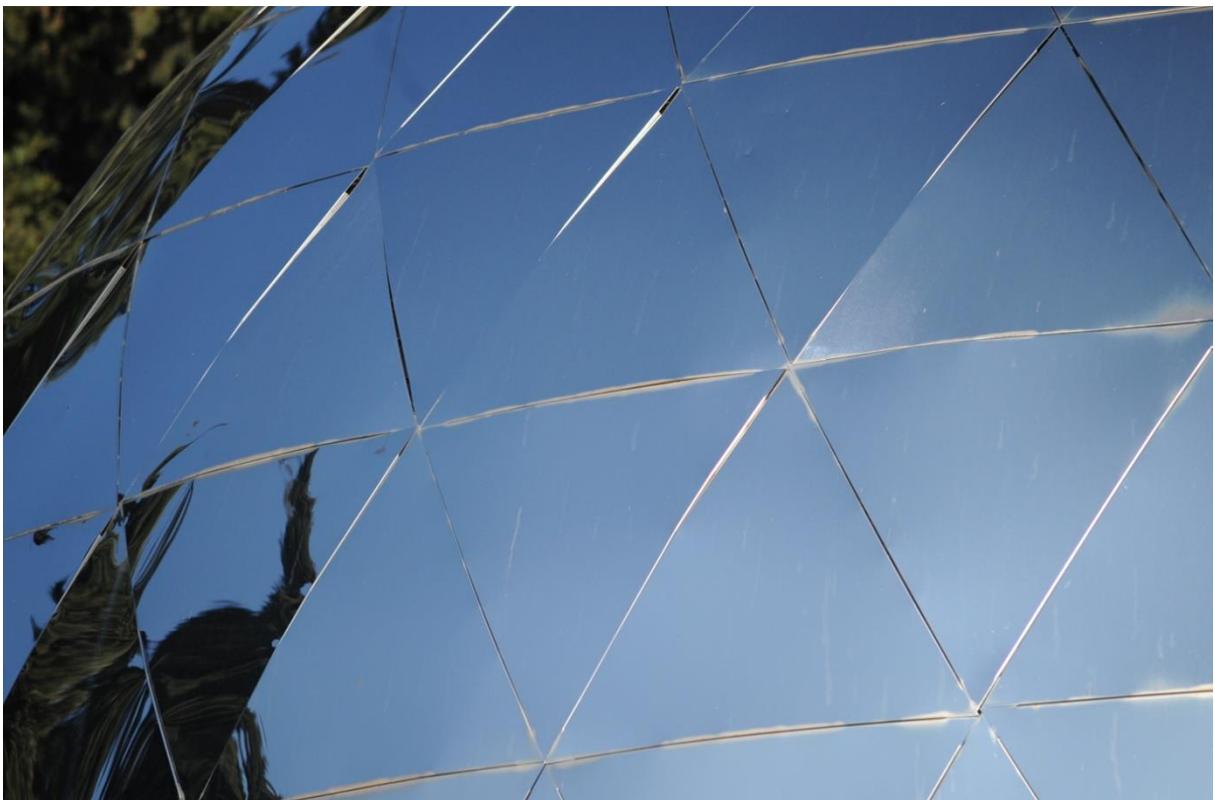

**Figure 3.** Detail of the hexagonal-only geodesic dome of the Latina's planetarium "Livio Gratton" (Italy).

**THE FSTF EDUCATIONAL PROJECT:**
Stimulated by five teachers (Department of Physical Education: M.G. Saviano, F. Truppa; Dept. of Mathematics and Physics: A. Curinga, R. Calvitti, the author), in march 2013 six Italian pupils aged 18-19 (G. Cardarelli, M. Garrido Cutiño, L. Gambacorta, I. Olivadese, M. Rosella, K. Tytler) explored the footballs' structure, familiarly known as *bucky-ball* [6, 13] still after *R. Buckminster Fuller* and his famous geodesic domes. According to the geometric construction of regular polygons, the students drew hexagons and pentagons on differently coloured plastic sheets (following the soccer balls' pattern: *white* hexagons and *black* pentagons) and then they cut them off. The cut polygons were then applied on a polystyrene ball with a diameter ($d = 25\ cm$) chosen a bit larger than the official football "size 5" ($d = 22 \div 23\ cm$) for practical purposes. The pupils used pins and black thread to emulate the seams of traditional footballs (Figs. 4, 5, 7 and 8). The students made two almost identical soccer-like balls, describing them in their State Exam (July 2013). As reported on the artefacts' book-shaped bases (Figs. 4, 5 and 8), this cross-sectional project involved knowledge and skills from four subjects: *Mathematics*, *Physical Education, Physics* and *Science*.

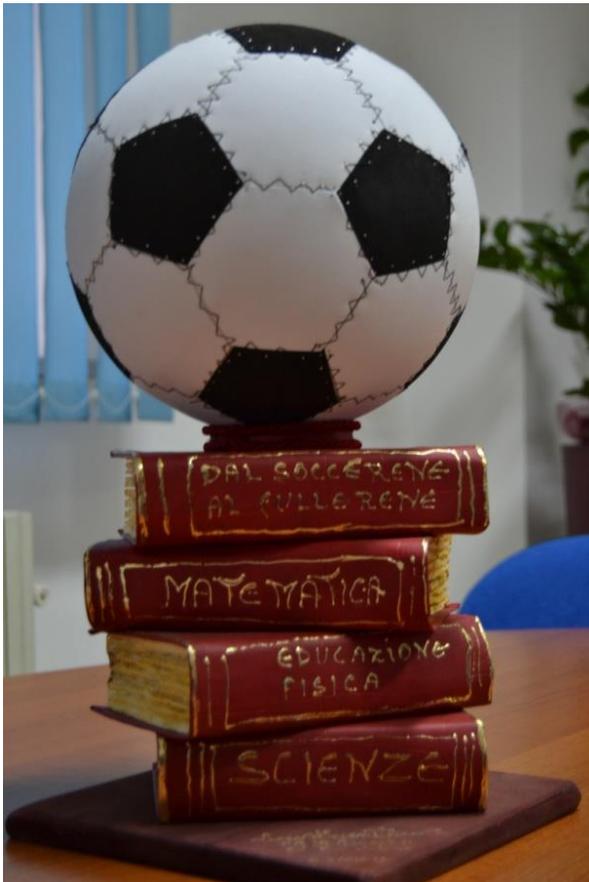 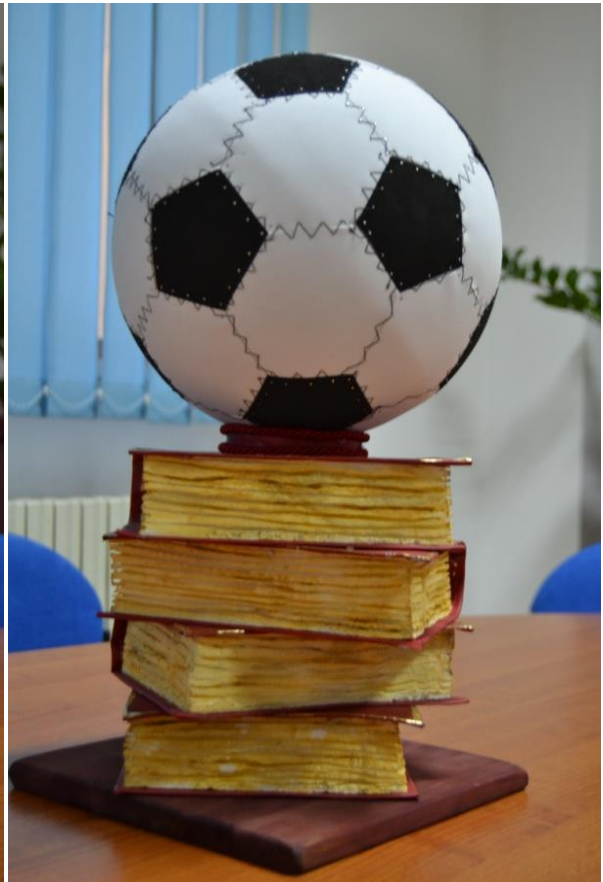

**Figure 4.** Front view of the 2013 FSTF artefact.      **Figure 5.** Back view of the 2013 FSTF artefact.

**THE MATHEMATICS OF THE FSTF:**
Once assimilated the geometry of the Carbon-60 molecule, the pupils were given two identical polystyrene spheres (each of a diameter of $25\ cm$), some plastic sheets (black and white), a $50\ m$ long black thread and around 2000 pins. Let us show the four computational steps autonomously taken by the pupils' team.

**Step 1. The unknown common edge**
The students considered that each of their two polystyrene spheres had the radius:

(1)   $r = 12.5\ cm$

and the spherical surface was:

(2)   $S_s = 4\pi r^2 = 4\pi \left(\frac{25}{2}\right)^2 = 625\pi \; cm^2$

The pupils realized that, since 20 hexagons and 12 pentagons must cover the whole spherical surface (Eq. (2)) they had to calculate the unknown common *x*-edge (or *x*-side, see Fig. 6) between the polygons.

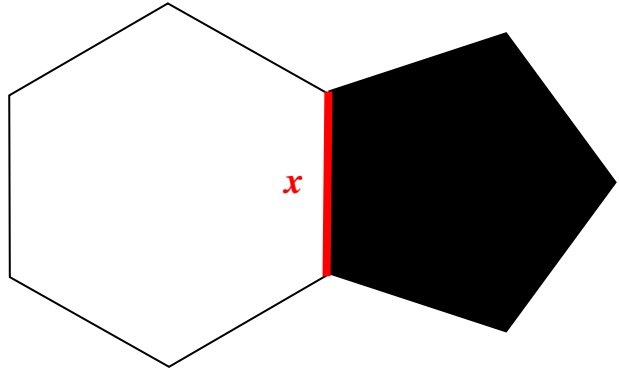

**Figure 6.** The common *x*-side of hexagons and pentagons.

### Step 2. The hexagons' side length
The students' team started from the hexagon's area, linked to the *x*-side by:

(3)   $A_h = \frac{3\sqrt{3}}{2} x^2$

The part of sphere covered by 20 hexagons was calculated as:

(4)   $S_h = 20 A_h = 30\sqrt{3} x^2$

### Step 3. The pentagons' side length
The students' team continued with the pentagon's area, linked to the *x*-side by:

(5)   $A_p = \frac{5}{4} x^2 \tan 54° = \frac{5}{4} x^2 \frac{1+\sqrt{5}}{\sqrt{10-2\sqrt{5}}}$

The part of sphere covered by 12 pentagons was calculated as:

(6)   $S_p = 12 A_p = 15 \frac{1+\sqrt{5}}{\sqrt{10-2\sqrt{5}}} x^2$

### Step 4. Value of the x-side
Since the sum of the hexagons' total surface (Eq. (4)) and of the pentagons' total surface (Eq. (6)) must cover the whole spherical surface (Eq. (2)), the pupils calculated the *x*-value as follows:

(7)   $S_s = S_h + S_p$

(8)   $625\pi = 30\sqrt{3} x^2 + 15 \frac{1+\sqrt{5}}{\sqrt{10-2\sqrt{5}}} x^2$

(9)   $x^2 \approx 27 \; cm^2$

(10)  $x \approx 5.2 \; cm$

The students were ready to cover their two polystyrene spheres by sewing together the white hexagons and black pentagons shown in Fig. 7, whose identical side was approximately 5.2 *cm* long.

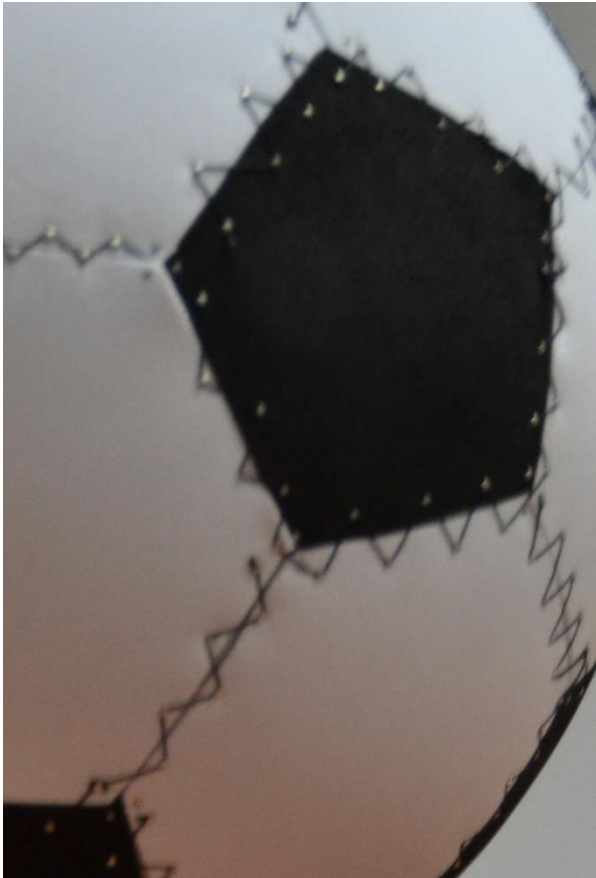 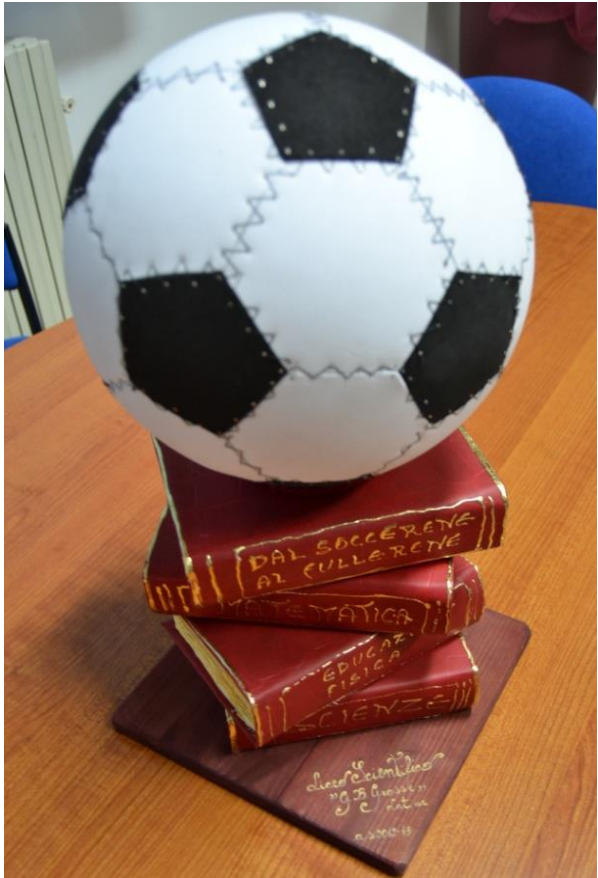

**Figure 7.** Detail of the 2013 FSTF artefact's sewings.   **Figure 8.** Top view of the 2013 FSTF artefact.

**CONCLUSIONS:**

Beyond the realization of a nice artefact (Fig. 8), the 2013 FSTF project's results were educational, cognitive, epistemological and heuristic. Among them, we recall:

(1) An innovative collaboration between traditionally distant Departments fixing an Italian school's weak point, i.e., the lack of enough interdisciplinary dialogue and activity.
(2) An application of the widely auspicated *learning by doing* [7] and *constructivism* [17] in Science [8], [10].
(3) The discussion of cross-sectional issues ([9], [11], [16], [18]) useful to improve the students' preparation for their state exam.
(4) A consideration on how even apparently easy sport questions, like the choice to play with a *ball*, can induce non-trivial reflections of physical order [14], [20].
(5) A path of gradual awareness about the difficulty to build, with flat surfaces, a convex solid approximating a sphere.
(6) The suggestion to extend the Italian pre-university curriculum to the *Archimedean solids*, besides the five classic (Platonic) polyhedra.

**ACKNOWLEDGMENTS:**

I'm grateful to the Organizing Committee of the IV ICMCS 2015, especially the Chairpersons Arjun K.Gupta and Krzysztof Drachal and the IJMSA Managing Editor Vinay Jha, for having given me the chance to present a pioneering Italian educational experience in the magnificent scenario of the AITCC in Thailand.
My thanks to the students Simone Ciotti and Valeria Volpe for their gift of a large amount of splendid pictures among which we have chosen the figures (respectively, the planetarium and the ball) in the paper.

**Mathematical Insights in the Pioneering Educational Project FSTF**

I also wish to thank both the FSTF Manager Francesco Truppa and the Grassi's Headmaster Loretta Tufo for their confidence in my ability of divulging the project's characteristics and achievements through papers and talks around the world.

**REFERENCES:**
1. M. Ausili, "The pontine project "From Soccerene to Fullerene" illustrated by Enzo Bonacci at the 100[th] SIF 2014 congress in Pisa," Latina in Vetrina (August 30, 2014). Via: http://www.latinainvetrina.it/content/il-progetto-pontino-"dal-soccerene-al-fullerene"-illustrato-dal-prof-enzo-bonacci-al-100°-co.
2. E. Bonacci, "The Planetarium experience in Latina (Italy)," In: Girep-Epec conference on Physics Education, University of Jyväskylä (August 1-5, 2011), JYFL Research Report No. 5, Physics Alive. Via: https://congress.cc.jyu.fi/girep2011/ preview/abstracts/long_abstract/1012_1320441491.pdf.
3. E. Bonacci, "The Planetarium experience in Latina," In: 97[th] National Congress of the Italian Physical Society, L'Aquila (September 26-30, 2011), Atticon6273 VI-C-1.
4. E. Bonacci, "The History of Science and Science Education: A Planetarium at School," In: Physics, Astronomy and Engineering. Critical Problems in the History of Science and Society, Šiauliai, The Scientia Socialis Press (2013) 141-146. Via: http://www.sisfa.org/wp-content/uploads/2013/12/SISFA2012Cont.pdf.
5. E. Bonacci, "The educational project FSTS in Latina," In: 100[th] conference of the Italian Physical Society, Atticon8187 VI-C-2, Pisa, Italy, September 22-26, 2014.
6. F.R.K. Chung, S. Sternberg, "Mathematics and the Buckyball," American Scientist 81(1) (1993) 57-61.
7. J. Dewey, *Democracy and education: An introduction to the philosophy of education*, New York, WLC Books, 2009 (original work published 1916).
8. F. Giuliano, "Rationalism and empiricism in action: the scientific method. Dogmatism or constructivism?," In: II workshop on Science Perception, Rome (April 19, 2013).
9. F. Giuliano, "Fullerenes: how the human imagination imitates Chemistry," Buongiorno Latina (January 16, 2014). Via: http://www.buongiornolatina.it/i-fullereni-ecco-come-la-fantasia-delluomo-imita-la-chimica-di-francesco-giuliano/.
10. F. Giuliano, "Razionalismo ed empirismo in azione: il metodo scientifico. Dogmatismo o costruttivismo?," In: AIF-Latina workshop on Science and Science Education, Latina, March 11, 2014.
11. F. Giuliano, "The synthesis of some chemical substances discovered by Italian Chemists," Etalia (July 24, 2014). Via: http://www.etalia.net/articles/a3e4613a9-60cd-4eb3-8e31-b2134d785113.
12. E. Knowles, J. Elliott, *The Oxford dictionary of new words*, Oxford, Oxford University Press, 1998.
13. O. Kuchment, "Some Short Histories of Common Scientific Terms: A Walk Through the Oxford English Dictionary," Science Editor 32(1) (2009) 23-24.
14. N. Ludwig, G. Guerrerio, *Science in Football: the secrets of Soccer revealed with Physics*, Bologna, Zanichelli, 2011. Via: http://online.scuola.zanichelli.it/chiavidilettura/files/2011/10/calcio_protetto.pdf.
15. J.S. Moore, "The Bucky Fuller Virtual Institute," BFVI (June 9, 1995). Via: http://www.u.arizona.edu/~shunter/bfvi.html.
16. J. Pedersen, "Nanotechnology Through History: Carbon-based Nanoparticles from Prehistory to Today," Sustainable-Nano (June 17, 2013). Via: http://sustainable-nano.com/2013/06/17/nanotechnology-through-history-carbon-based-nanoparticles-from-prehistory-to-today/.
17. J. Piaget, *The construction of reality in the child*, New York, Basic Books, 1954.
18. E.F. Porter, "A Whole new Ball Game," St Louis Post-Dispatch (November 8, 1993). Via: http://www.questia.com/ newspaper/1P2-32837795/a-whole-new-ball-game.
19. D. Shen, *Approaches to Soccerene ($I_h C60$) and Other Carbon Spheres*, Los Angeles, University of California, 1990.
20. L.R.G. Treloar, *The Physics of Rubber Elasticity*, Oxford, Oxford University Press, 1975.
21. E.W. Weisstein, "Polyhedral Formula," From MathWorld - A Wolfram Web Resource. Via: http://mathworld.wolfram.com/PolyhedralFormula.html.